**Regional Boundary Asymptotic Gradient Full-Order Observer in Distributed Parabolic Systems**

Zinah A. Khalid[1] and Raheam A. Al-Saphory[2,]

[1,2]Department of Mathematics, College of Education for Pure Sciences, Tikrit University, Tikrit, Iraq.

saphory@hotmail.com[2]  zeena_assif@hotmail.com[1]

**Abstract**

The purpose of this paper is to explore the concept of the regional boundary asymptotic gradient full order observer (RBAGFO-observer) in connection with the characterizations of sensors structures. Then, we present various results related to different types of measurements, domains and boundary conditions for distributed parameter systems (DPS$_S$) in parabolic systems problem.  The considered approach of this work is derived from Luenberger observer theory which is enable to estimate asymptotically the state gradient of the original system on a sub-region of the domain boundary $\partial\Omega$  in order that the RBAGFO-observability notion to be achieved. We also show that there exists a dynamical system for the considered system is not BAGFO-observer in the usual sense, but it may be regional RBAGFO-observer.

**1. Introduction**

There are many situations in modern technology in which it is necessary to estimate the state of a dynamic system using only the measured input and output data of the system [1]. An observer is a dynamic system Ŝ which to estimates the state of another considered system S using only the measured input and output letter. If the order of Ŝ is equal to the order of S, the observer is called full-order state observer [2-4].

The asymptotic observer theory explored by Luenberger in [5] for finite-dimensional linear systems and extended infinite-dimensional distributed parameter systems govern by strongly continuous semi-group in Hibert space by Gressang and Lamont as in [6]. The study of this approach via another variable like sensors and actuators developed by El-Jai *et al*. as in ref.s [2-7] in order to achieve asymptotic observability.

One of the most important approach in system theory is focused on the reconstruction of the state of the system from the knowledge of dynamic system and the output function on a sub-region $\omega$ of a spatial domain $\Omega$ this problem is called regional observability problem has been received much attention as in [8-10].

An extension of this notion has been given in [11-12] to the regional gradient case. The regional asymptotic notion has been introduced and developed by Al-Saphory and El- Jai in [13-14]. Thus, this notion consists in studying the asymptotic behavior of the system in an internal sub-region $\omega$ of a spatial domain $\Omega$.

Thus, the asymptotic regional state reconstruction studied and developed in [15-16]  and extended to the regional asymptotic gradient full-order observer (RGFO-observer) which allows estimating the state gradient of the original system.

The purpose of this paper is to study and examine the concept of RBAGFO-observer by using the choice of sensors. The principle reason for considering this case is that, in the first time, the existent of a dynamical system which is observed asymptotically the gradient of the system state on some boundary region Γ ⊂ $\partial\Omega$ [17]. In second time, it is closer to a real situation, the treatment of water by using a bioreactor where the objective is to observe the concentration of substrate at the boundary output of the bioreactor in order the water regulation is achieved (for example see figure 1).

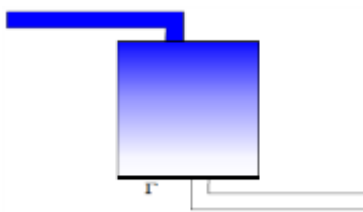

**Fig. 1:** Observation of substrate concentration at the boundary output of the reactor.

The outline of this paper is organized as follows: Section 2 is devoted to the problem statement and some basic concept related to the regional boundary gradient stability (RBAG-stability), regional boundary asymptotic gradient detectability (RBAG-detectability). Section 3, we focus on RBAGFO-observer so we introduced and characterize the existing of RBAGFO-observer to provide RBAGFO-estimator of gradient state for the original system in terms of sensors structure. In the last section, we have been applied these result to the two (DPS$_S$) for different zone and pointwise sensors case.

## 2. Problem Formulation and Preliminaries

This section present considered system and formulation of the problem with some definitions and characterizations, which is related to the present work.

### 2.1 Problem Statement

Let $\Omega$ be an open bounded subset of $R^n$ with boundary $\partial\Omega$ and $\Gamma$ be a region subset of $\partial\Omega$. We denote $Q = \Omega \times ]0,\infty[$ and $\Sigma = \partial\Omega \times ]0,\infty[$. Consider the parabolic system which is described by the following state-space equation

$$\begin{cases} \frac{\partial x}{\partial t}(\xi,t) = Ax(\xi,t) + Bu(t) & Q \\ x(\xi,0) = x_0(\mu) & \bar{\Omega} \\ \frac{\partial x}{\partial v}(\eta,t) = 0 & \Sigma \end{cases} \qquad (1)$$

augmented with the output function

$$y(.,t) = Cx(.,t) \qquad (2)$$

• The separable Hilbert spaces are $X$, $U$ and $\mathcal{O}$ where $X = H^1(\bar{\Omega})$ is the state space, $U = L^2(0,T,R^p)$ is the control space and $\mathcal{O} = L^2(0,T,R^q)$ is the observation space, where $p$ and $q$ are the numbers of actuators and sensors.

• $A = \sum_{i,j=1}^n \frac{\partial}{\partial x_j}(a_{ij} \frac{\partial}{\partial x_j})$ with $a_{ij} \in D(\bar{A})$ (the domain of $\bar{A}$) is a second-order linear differential operator, which generates a strongly continuous semi-group $(S_A(t))_{t\geq 0}$ on the state space $X$ and is self-adjoint with compact resolvent [18].

• The operators $B \in L(R^p, X)$ and $C \in L(X, R^q)$ depend on the structure of actuators and sensors as in [19] **(**figure 2**)**

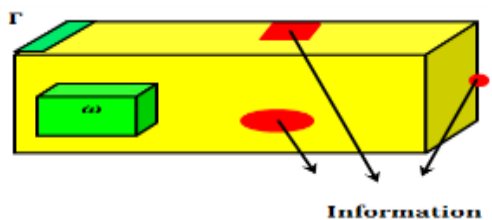

**Fig. 2:** Mathematical model: domain Ω, region Γ, and sensors Locations.

- The mathematical model in figure 2 is more general and complex than the the spatial case of real model in figure 1.

- Under the given assumptions above, the system (1) has a unique solution given by the following form [20-21].

$$x(\xi,t) = S_A(t)x_0(\xi) + \int_0^t S_A(t-s)Bu(s)\,ds \tag{3}$$

- The problem is how to reconstruct a dynamical system may be called estimator for the current state gradient in a given region on Γ, and to give a sufficient condition for the existence of a RBAGFO-observer.

- The initial state $x_0$ and its gradient $\nabla x_0$ are supposed to be unknown, the problem concerns the reconstruction of the initial gradient $\nabla x_0$ on the region Γ of the system domain $\partial\Omega$.

- Now, we consider the operator $K$ given by the form

$$K: X \to \mathcal{O} \tag{4}$$

$$x \to CS_A(.)x$$

where $K$ is a bounded linear operator as in [7, 22, 24]. And the adjoint operator $K^*$ of $K$ is defined by

$K^*: \mathcal{O} \to X$, and represented by the form

$$K^*y^* = \int_0^t S_A^*(s)C^*y^*(s)ds \tag{5}$$

- The operator $\nabla$ denotes the gradient is given by

$$\begin{cases} \nabla: H^1(\Omega) \to (H^1(\Omega))^n \\ x \to \nabla_x = \left(\frac{\partial x}{\partial \xi_1}, \ldots, \frac{\partial x}{\partial \xi_n}\right) \end{cases} \tag{6}$$

and, the adjoint of $\nabla$ denotes by $\nabla^*$ is given by

$$\begin{cases} \nabla^*: (H^1(\Omega))^n \to H^1(\Omega) \\ x \to \nabla_x^* = v \end{cases} \tag{7}$$

where $v$ is a solution of the Dirichlet problem

$$\begin{cases} \Delta_v = -div(x) & \Omega \\ v = 0 & \partial\Omega \end{cases} \tag{8}$$

- The trace operator of order zero is described by [23]

$$\gamma_0: H^1(\Omega) \to H^{1/2}(\partial\Omega) \tag{9}$$

which is a linear, subjective, and continuous [20]. Thus, the extension of the trace operator [23] which is denoted by $\gamma$ defined as

$$\gamma: (H^1(\Omega))^n \to (H^{1/2}(\partial\Omega))^n \tag{10}$$

and the adjoints are respectively given by $\gamma_0^*, \gamma^*$.

- For a region Γ of $\partial\Omega$, we define the gradient restriction operator by the form

$$\chi_\Gamma: (H^{1/2}(\partial\Omega))^n \to (H^{1/2}(\Gamma))^n \qquad (11)$$

where the adjoint of $\chi_\Gamma$ denotes by $\chi_\Gamma^*$ is defined by

$$\chi_\Gamma^*: (H^{1/2}(\Gamma))^n \to (H^{1/2}(\partial\Omega))^n \qquad (12)$$

• Finally, we denote the operator $H_{\Gamma^*G} = \chi_\Gamma \gamma \nabla K^*$ from $\mathcal{O}$ into $(H^{1/2}(\Gamma))^n$ and the adjoint of this operator given by $H_{\Gamma^*G}^* = K\nabla^*\gamma^*\chi_\Gamma^*$.

• For a sub-region $\Gamma$ of the boundary regional boundary gradient $\Gamma_{RBG} = \Gamma^*$.

Now, the problem is how to build an approach which observe (estimates) regional state gradient on a region $\Gamma^*$ of the boundary $\partial\Omega \subset \bar{\Omega}$ asymptotically by using a dynamic system (an observer) in full order case only may be called full-order observer in region $\Gamma^*$. The important of an observer is that to estimates all the gradient of state variables, regardless of whether some are available for direct measurements or not [1].

## 2.2 $\Gamma^*$G-observability and $\Gamma^*$ AG-detectability

This sub section devotes the relation between the concept of RBG-observability and RBAG-detectability on $\Gamma^*$. As well known the observability [19-22] and asymptotic observability [1-5, 7, 24] are important concepts to estimate the unknown state of the considered dynamic system from the input and output functions. Thus, These notions are studied and introduced to the DPSs analysis with different characterizations by El-Jai, Zerrik and Al-Saphory *et al.* in many paper for example [8-18, 28-31] in connection with strategic sensors.

• The systems (1)-(2) are said to be exactly regionally boundary gradient observable (E$\Gamma^*$G-observable) on $\Gamma^*$ if

$$Im\, H = Im\chi_{\Gamma^*}\nabla K^* = (H^{1/2}(\Gamma^*))^n$$

• The systems (1)-(2) are said to be weakly regionally boundary gradient observable (W$\Gamma^*$G-observable) on $\Gamma^*$ if

$$\overline{Im\, H} = \overline{Im\chi_{\Gamma^*}\nabla K^*} = (H^{1/2}(\Gamma^*))^n$$

It is equivalent to say that the systems (1)-(2) are W$\Gamma^*$G-observable if

$$Ker\, H^* = ker K\, \nabla^*\chi_{\Gamma^*} = \{0\}$$

• If the systems (1)-(2) are is W$\Gamma^*$G-observable, then $x_0(\xi, 0)$ is given by

$$x_0 = (K^*K)^{-1}K^*y = K^\dagger y,$$

where $K^\dagger$ is the pseudo-inverse of the operator $K$ [9-10].

• A sensor $(D, f)$ is a regional boundary gradient strategic ($\Gamma^*$G-strategic) on $\Gamma^*$ if the observed system is W$\Gamma^*$G-observable.

• The measurements can be obtained by the use of zone or pointwise sensors, which may be located in $\Omega$ or $\partial\Omega$ [24].

• Then, according to the choice of the parameters $D_i$ and $f_i$, we have different types of sensors:

• It may be zone, if $D_i \subset \bar{\Omega}$ and $f_i \in L^2(D_i)$. In this case, the operator $C$ is bounded, and the output function (2) may be given by the form

$$y(t) = \int_{D_i} z(\xi,t) f_i(\xi) \, d\xi \tag{13}$$

- It may be pointwise, if $D_i = \{b_i\}$ with $b_i \in \bar{\Omega}$ and $f = \delta(.-b_i)$ where $\delta$ is a Dirac mass concentrated in $b$ [14, 21, 28]. In this case, the operator $C$ is un bounded, and the output function (2) may be given by the form

$$y(t) = \int_\Omega x(\xi,t) \delta_{b_i}(\xi - b_i) \, d\xi \tag{14}$$

- It may be boundary zone, if $\Gamma_i \subset \partial\Omega$ and $f_i \in L^2(\Gamma_i)$, the output function (2) may be given by the form

$$y(t) = \int_{\Gamma_i} x(\eta,t) f_i(\eta) \, d\eta \tag{15}$$

**Definition 2.1:** The semi-group $(S_A(t))_{t \geq 0}$ is regionally boundary asymptotically gradient stable ($\Gamma^*$AG-stable) on $\Gamma^*$, if and only if for some positive constants $M_{\Gamma^*}$, $\alpha_{\Gamma^*}$, we have

$$\|\chi_{\Gamma^*} \gamma \nabla S_A(.)\|_{L((H^{1/2}(\Gamma^*))^n, H^1(\bar{\Omega}))} \leq M_{\Gamma^*} e^{\alpha_{\Gamma^*}}, \forall t \geq 0.$$

**Remark 2.2:** If the semi-group $(S_A(t))_{t \geq 0}$ is $\Gamma^*$AG-stable on $(H^{1/2}(\Gamma^*))^n$, then for all $x_o \in H^1(\Omega)$, the solution of the associated system satisfies

$$\lim_{t \to \infty} \|\chi_{\Gamma^*} \gamma \nabla x(.,t)\|_{(H^{1/2}(\Gamma^*))^n} = \lim_{t \to \infty} \|\chi_{\Gamma^*} \gamma \nabla S_A(t) x_0\|_{(H^{1/2}(\Gamma^*))^n}$$

$$= 0 \tag{16}$$

**Definition 2.3:** The system (1) is said to be $\Gamma^*$AG-stable on $\Gamma^*$ if the operator $A$ generates a semi-group which is $\Gamma^*$AG-stable on the space $(H^{1/2}(\Gamma^*))^n$.

**Definition 2.4:** The system (1)-(2) is said to be regionally boundary asymptotically gradient detectable ($\Gamma^*$AG-detectable) on $\Gamma^*$, if there exists an operator $H_{\Gamma^*AG}: R^q \to (H^{1/2}(\Gamma^*))^n$, such that the operator $(A - H_{\Gamma^*AG} C)$ generates a strongly continuous semi-group $(S_{H_{\Gamma^*AG}}(t))_{t \geq 0}$, which is $\Gamma^*$AG-stable on $(H^{1/2}(\Gamma^*))^n$.

**Proposition 2.5:** If the system (1)-(2) is $\Gamma^*$G-observable on $\Gamma^*$, then it is $\Gamma^*$AG-detectable on $\Gamma^*$. This results gives the following inequality: $\exists k_{\Gamma^*AG} > 0$, such that

$$\|\chi_{\Gamma^*} \gamma \nabla S_A(.)x\|_{(H^{1/2}(\Gamma^*))^n} \leq k_{\Gamma^*AG} \|C S_A(.)x\|_{L^2(0,\infty,O)}, \tag{17}$$

for all $z \in (H^{1/2}(\Gamma^*))^n$.

**Proof:** We conclude the proof of this proposition is conclude from the results on observability considering $\chi_{\Gamma^*} \nabla K^*$. We have the following forms [21, 24]

1. $Im f \subset Im g$.

2. There exists $k > 0$, such that

$$\|f^* x^*\|_{E^*} \leq k \|g^* x^*\|_{F^*}, \text{for all } x^* \in G^*$$

From the right hand said of above inequality $k_{\Gamma^*AG} \|g^* x^*\|_{F^*}$, there exists $M_{\Gamma^*AG}, \omega_{\Gamma^*AG} > 0$ with $k_{\Gamma^*AG} < M_{\Gamma^*AG}$, such that

$$k_{\Gamma^*AG} \|g^* x\|_{F^*} \leq M_{\Gamma^*AG} e^{-\omega_{\Gamma^*AG} t} \|x^*\|_{F^*}$$

where $E, F$ and $G$ be a reflexive Banach spaces and $f \in L(E,G)$, $g \in L(F,G)$. If we apply this result, considered

$$E = G = (H^{1/2}(\Gamma^*))^n, \ F = \mathcal{O}, \ f = Id_{(H^{1/2}(\Gamma^*))^n}$$

and

$$g = S_A^*(.)\chi_\Gamma^* \gamma^* \nabla^* C^*$$

where $S_A(.)$ is a strongly continuous semi-group generates by $A$, which is $\Gamma^*$AG-stable on $\Gamma^*$, then it is $\Gamma^*$AG-detectable on $\Gamma^*$. Thus, the notion of $\Gamma^*$AG-detectability is a weaker property than the $\Gamma^*$G-observability [25-26].

## 3. Sensors and $\Gamma^*$AGFO-Observer

In this section we present the sufficient conditions which are guarantee the existence of regional boundary asymptotic gradient full order observer ($\Gamma^*$AGFO-Observer) on $\Gamma^*$ which allows to construct a $\Gamma^*$AGFO-estimator on $\Gamma^*$ of the state $\chi_{\Gamma^*}\gamma\nabla Tx(\xi,t)$.

### 3.2 Definitions and characterizations

**Definition 3.1:** Suppose there exists a dynamical system with state $z(.,t) \in Z$ given by

$$\begin{cases} \frac{\partial \hat{z}}{\partial t}(\xi,t) = A\hat{z}(\xi,t) + Bu(t) + H_{\Gamma^*AG}\big(Cx(\xi,t) - C\hat{z}(\xi,t)\big) & Q \\ \hat{z}(\xi,0) = \hat{z}_o(\xi) & \bar{\Omega} \\ \hat{z}(\eta,t) = 0 & \Sigma \end{cases} \quad (18)$$

In this case the operator $F_{\Gamma^*AG}$ in general case [13] is given by $F_{\Gamma^*AG} = A - H_{\Gamma^*AG}C$ where $T = I$ the identity operator. Thus the operator $A - H_{\Gamma^*AG}C$ generate a strongly continuous semi-group $(S_{A-H_{\Gamma^*AG}C}(t))_{t\geq 0}$ on separable Hilbert space $Z$ which is $\Gamma^*$AG-stable.

Thus, $\exists \ M_{A-H_{\Gamma^*AG}C}, \alpha_{A-H_{\Gamma^*AG}C} > 0$ such that

$$\left\|S_{A-H_{\Gamma^*AG}C}(.)\right\| \leq M_{A-H_{\Gamma^*AG}C} e^{-\alpha_{A-H_{\Gamma^*AG}C}t}, \ \forall t \geq 0.$$

and let $G_{\Gamma^*AG} \in L(U,Z), H_{\Gamma^*AG} \in L(\mathcal{O},Z)$ such that solution of (18) similar to (3)

$$z(\xi,t) = S_{A-H_{\Gamma^*AG}C}(t)z(\xi) + \left[\int_0^t S_{A-H_{\Gamma^*AG}C}(t-\tau) \, Bu(\tau) \, H_{\Gamma^*AG} y(\tau)\right] d\tau.$$

**Definition 3.2:** The system (18) defines $\Gamma^*$AGFO-estimator such that

$$z(\xi,t) = \chi_\Gamma \nabla Tx(\xi,t) = Ix(\xi,t) \in (H^{1/2}(\Gamma^*))^n$$

where $x(\xi,t)$ is the solution of the systems (1)-(2) if

$$\lim_{t\to\infty} \|z(.,t) - x(\xi,t)\|_{(H^{1/2}(\Gamma^*))^n} = 0,$$

and $\chi_\Gamma \nabla I$ maps $D(A)$ into $D(A - H_{\Gamma^*AG}C)$ where $z(\xi,t)$ is the solution of system (18).

**Remark 3.3:** The dynamic system (18) specifies $\Gamma^*$AGFO-observer of the systems given by (1)-(2) if the following holds:

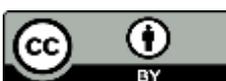

1-There exists $M_{\Gamma^*AGFO} \in L\left(R, \left(H^{1/2}(\Gamma^*)\right)^n\right)$ and $N_{\Gamma^*AGFO} \in L\left(\left(H^{1/2}(\Gamma^*)\right)^n\right)$ such that

  i- $M_{\Gamma^*AGFO}C + N_{\Gamma^*AGFO} = I_{\Gamma^*AGFO}$

  ii- $A - F_{\Gamma^*AGFO} = H_{\Gamma^*AGFO}C$ and $G_{\Gamma^*AGFO} = B$.

3- The system (18) defines $\Gamma^*$AGFO-estimator for $x(\xi,t)$.

The object of $\Gamma^*$AGFO-observer is to provide an approximation to the original system state gradient. This approximation is given by

$$\hat{x}(t) = M_{\Gamma^*AGFO}y(t) + N_{\Gamma^*AGFO}z(t).$$

**Definition 3.4:** The systems (1)-(2) are regionally boundary asymptotic gradient full order observable ($\Gamma^*$AGFO-observable) on $\Gamma^*$, if there exists a dynamic system which is $\Gamma^*$AGFO-observer for the considered system.

### 3.2 $\Gamma^*$AGFO-Observer reconstruction method

In this case, we need to consider $\chi_{\Gamma^*AGFO}\nabla T = I_{\Gamma^*AGFO}$ and $Z = X$, then the operator observer equation becomes as $F_{\Gamma^*AGFO} = A - H_{\Gamma^*AGFO}C$ where $A$ and $C$ are known. Thus, the operator $H_{\Gamma^*AGFO}$ must be determined such that the operator $F_{\Gamma^*AGFO}$ is $\Gamma^*$AG-stable. This observer is an extension of asymptotic observer as in [14-16, 24-31]. Now consider again system (1) together with output function (2) described by the following form

$$\begin{cases} \frac{\partial x}{\partial t} = Ax(\xi,t) + Bu(t) & Q \\ x(\xi,0) = x_0(\xi) & \bar{\Omega} \\ x(\eta,t) = 0 & \Sigma \\ y(t) = Cx(.,t) & Q \end{cases} \quad (19)$$

Let $\Gamma^*$ be a given sub-domain of $\bar{\Omega}$ and suppose that $I_{\Gamma^*AGFO} \in \mathcal{L}(H^1(\Omega))^n$, and $\chi_{\Gamma^*AGFO}\nabla Tx(\xi,t) = \chi_{\Gamma^*AGFO}\nabla x(\xi,t)$ there exists a system with state $z(\xi,t)$ such that

$$z(\xi,t) = \chi_{\Gamma^*AGFO}\nabla Tx(\xi,t) = T_{\Gamma^*AGFO}x(\xi,t) \qquad (20)$$

(with $T_{\Gamma^*AGFO} = I_{\Gamma^*AGFO}$ where $I_{\Gamma^*AG}$ is the identity operator with respect to $\Gamma^*$AGFO-estimator. Then

$$z(\xi,t) = I_{\Gamma^*AGFO}x(\xi,t) = x(\xi,t) \qquad (21)$$

From equation (2) and (21) we have

$$\begin{bmatrix} y \\ z \end{bmatrix} = \begin{bmatrix} C \\ I_{\Gamma^*AGFO} \end{bmatrix} x$$

If we assume that there exist two bounded linear operators

$$M_{\Gamma^*AGFO}: \vartheta \to ((H^{1/2}(\Gamma^*))^n$$

and

$$N_{\Gamma^*AGFO}: ((H^{1/2}(\Gamma^*))^n \to ((H^{1/2}(\Gamma^*))^n$$

such that

$$M_{\Gamma^*AGFO}C + N_{\Gamma^*AGFO}T_{\Gamma^*AGFO} = I_{\Gamma^*AGFO}$$

then by deriving $z(\xi,t)$ in (20) we have

$$\frac{\partial z}{\partial t}(\xi,t) = I_{\Gamma^*AGFO}\frac{\partial x}{\partial t}(\xi,t) = \chi_{\Gamma^*AGFO}\nabla T A x(\xi,t) + \chi_{\Gamma^*AGFO}\nabla T B u(t)$$

$$= \chi_{\Gamma^*}\nabla I_{\Gamma^*AGFO} A\, M_{\Gamma^*AGFO} y(\xi,t) + \chi_{\Gamma^*}\nabla I_{\Gamma^*AGFO} A N_{\Gamma^*AGFO} z(\xi,t)$$

$$+ \chi_{\Gamma^*}\nabla I_{\Gamma^*AGFO} B u(t)$$

Since the operator $T_{\Gamma^*AGFO} = I_{\Gamma^*AGFO}$, then we have

$$\frac{\partial z}{\partial t}(\xi,t) = A N_{\Gamma^*AGFO}\, z(\xi,t) + B u(t) + A M_{\Gamma^*AGFO} y(\xi,t)$$

Therefore

$$\frac{\partial \hat{z}}{\partial t}(\xi,t) = F_{\Gamma^*G}\hat{z}(\xi,t) + G_{\Gamma^*G} u(t) + H_{\Gamma^*G}\, y(.,t)$$

and since $A - F_{\Gamma^*AGFO} = H_{\Gamma^*AGFO} C$ and $G_{\Gamma^*AGFO} = B$ then we have

$$\begin{cases} \frac{\partial \hat{z}}{\partial t}(\xi,t) = A\hat{z}(\xi,t) + Bu(t) + H_{\Gamma^*AGFO}(Cx(\xi,t) - cC(\xi,t)) & Q \\ \hat{z}(\xi,0) = \hat{z}_0(\xi) & \overline{\Omega} \\ \hat{z}(\eta,t) = 0 & \Sigma \end{cases} \quad (22)$$

Let us consider a complete sets of eigenfunctions $\varphi_{nj}$ in $(H^1(\Omega))^n$ orthonormal to $(H^{1/2}(\Gamma^*))^n$ associated with the eigenvalue $\lambda_n$ of multiplicity $r_n$ and suppose the system (1) has unstable mode. Then, the sufficient condition of an $\Gamma^*$AGFO-observer is formulated in the following main result.

**Theorem 3.5:** Suppose that there are $q$ zone sensors $(D_i, f_i)_{1\le i\le q}$ and the spectrum of $A$ contains $J$ eigenvalues with non-negative real parts. Then the dynamic system (22) is $\Gamma^*$AGFO-observer on $\Gamma^*$ for the system (19), that is $\lim_{t\to\infty}[z(\xi,t) - \hat{z}(\xi,t)] = 0$, if :

1-There exists

$M_{\Gamma^*AGFO} \in L(R^q, (H^{1/2}(\Gamma^*))^n)$ and $N_{\Gamma^*AGFO} \in L((H^{1/2}(\Gamma^*))^n)$ such that

$\quad M_{\Gamma^*AGFO} C + N_{\Gamma^*AGFO} = I_{\Gamma^*AGFO}$.

2- $A - F_{\Gamma^*AGFO} = H_{\Gamma^*AGFO} C$, $G_{\Gamma^*AGFO} = B$.

3- $q \ge m$

4- rank $G_m = m_m$, $\forall m, m = 1,\dots, J$ with

$$G_m = (G_m)_{ij} = \begin{cases} \langle \psi_{mj}(b_i), f_i(.) \rangle_{L^2(Di)} \\ \psi_{mj}(b_i) \\ \langle \frac{\partial \psi_{mj}}{\partial v}, f_i(.) \rangle L^2(\Gamma^*_i) \end{cases}$$

where sup $m_m = m < \infty$ and $j = 1,\dots, m_m$.

**Proof:**

**First step:** The proof is limited to the case of pointwise sensors. Under the assumptions of section 2, the system (1) can be decomposed by the projections $P$ and $I - P$ on two parts, unstable and stable. The state vector may be given by $x(\xi,t) = [x_1(\xi,t), x_2(\xi,t)]^t$ where $x_1(\xi,t)$ is the state component of the unstable part of the system (1) may be written in the form

$$\begin{cases} \frac{\partial x_1}{\partial t}(\xi,t) = A_1 x_1(\xi,t) + PBu(t) & Q \\ x_1(\xi,0) = x_{01}(\xi) & \overline{\Omega} \\ x_1(\eta,t) = 0 & \Theta \end{cases} \quad (23)$$

and $x_2(\xi,t)$ is the component state of the stable part of the system (1) given by

$$\begin{cases} \frac{\partial x_2}{\partial t}(\xi,t) = A_2 x_2(\xi,t) + (I-P)Bu(t) & Q \\ x_2(\xi,0) = x_{02}(\xi) & \overline{\Omega} \\ x_2(\eta,t) = 0 & \Theta \end{cases} \quad (24)$$

The operator $A_1$ is represented by a matrix of order $(\sum_{m=1}^J m_m, \sum_{m=1}^J m_m)$ given by $A_1 = \text{diag}[\lambda^1, \ldots, \lambda^1, \lambda^2, \ldots, \lambda^2, \ldots, \lambda_J, \ldots, \lambda_J]$ and $PB = [G_1^{t_r}, G_2^{t_r}, \ldots, G_J^{t_r}]$. The condition (4) of this theorem, allows that the suit $(D_i, f_i)_{1 \le i \le q}$ of sensors is $\Gamma^*G$-strategic for the unstable part of the system (1), the subsystem (23) is $W\Gamma^*G$ - observable [11] and since it is finite dimensional, then it is $E\Gamma^*G$ observable [27]. Therefore it is asymptotically $\Gamma^*AG$-detectable, and hence there exists an operator $H^1_{\Gamma^*AG}$ such that $(A_1 - H^1_{\Gamma^*AG} C)$ which is satisfied the following:

$\exists, M_{A1-H^1_{\Gamma^*AG}C}, \alpha_{A1-H^1_{\Gamma^*AG}C} > 0$ such that

$$\left\| e^{(A1-H^1_{\Gamma^*AG}C)t} \right\|_{(H^1(\Gamma^*))^n} \le M_{A1-H^1_{\Gamma^*AG}C} e^{-\alpha_{A1-H^1_{\Gamma^*AG}C}t}$$

and then we have

$$\|x_1(.,t)\|_{(H^1(\Gamma^*))^n} \le M_{A1-H^1_{\Gamma^*AG}C} e^{-\alpha_{A1-H^1_{\Gamma^*}C}} \|Px_0(.)\|_{(H^1(\Gamma^*))^n}$$

Since the semi-group generated by the operator $A_2$ is stable on $(H^{1/2}(\Gamma^*))^n$, then there exist $M_{A2-H^2_{\Gamma^*AG}C}, \alpha_{A2-H^2_{\Gamma^*AG}C} > 0$ [5] such that

$$\|x_2(.,t)\|_{(H^1(\Gamma^*))^n} \le M_{A1-H^1_{\Gamma^*AG}C} e^{-\alpha_{A2-H^2_{\Gamma^*AG}C(t)}} \|(I-P)\|_{(H^{1/2}(\Gamma^*))^n}$$

$$+ \int_0^t M_{M_{A1-H^1_{\Gamma^*AG}C}}, e^{-\alpha_{A2-H^2_{\Gamma^*AG}C(t)}} H^2_{\Gamma^*AG} \ \|(I-P)x_0(.)\|_{(H^{1/2}(\Gamma^*))^n} \|u(\tau)\| d\tau$$

and therefor $x(\xi,t) \to 0$ when $t \to \infty$. Finally, the system (23) are asymptotically $\Gamma^*AG$ -detectable.

**Second step:** From equation (21), we have $z(\xi,t) = x(\xi,t)$ with the observer error is given by the following form

$$e(\xi,t) = z(\xi,t) - \hat{z}(\xi,t)$$

where $\hat{z}(\xi,t)$ is a solution of the dynamic system (22). Derive the above equation, and by using equation (19) and condition 2, we can get the following forms

$$\frac{\partial \varepsilon}{\partial t}(\xi,t) = \frac{\partial z}{\partial t}(\xi,t) - \frac{\partial \hat{z}}{\partial t}(\xi,t) = \frac{\partial x}{\partial t}(\xi,t) - \frac{\partial \hat{z}}{\partial t}(\xi,t)$$

$$= Ax(\xi,t) + Bu(t) - F_{\Gamma^*AGFO} \ \hat{z}(\xi,t)$$

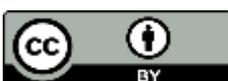

$$- G_{\Gamma^*} u(t) - H_{\Gamma^*\text{AGFO}} Cx(\xi,t)$$

$$= Ax(\xi,t) - (A - H_{\Gamma^*\text{AGFO}G} C)\hat{z}(\xi,t) - H_{\Gamma^*\text{AGFO}} Cx(\xi,t)$$

$$= (A - H_{\Gamma^*\text{AGFO}} C)(z(\xi,t) - \hat{z}(\xi,t))$$

$$= (A - H_{\Gamma^*\text{AGFO}} C) e(\xi,t)$$

Thus, from the first part of this proof we obtain

$$e(\xi,t) = (A - H_{\Gamma^*\text{AGFO}} C) e(0,t)$$

is asymptotically $\Gamma^*G$-stable with

$$e(0,t) = z_0(\xi) - \hat{z_0}(\xi)$$

Then we have

$$\|e(\xi,t)\|_{(H^{1/2}(\Gamma^*))^\eta} \leq M_{A-H_{\Gamma^*\text{AGFO}}C} e^{-\alpha A - H_{\Gamma^*\text{AGFO}}C^t}.$$

$$\|z_0(\xi) - \hat{z_0}(\xi)\|_{(H^{1/2}(\Gamma^*))^\eta}$$

therefore $\lim_{t \to \infty} e(\xi,t)$. Now, let the approximate solution to the gradient state of the original system is

$$\hat{x}(\xi,t) = M_{\Gamma^*\text{AGFO}} y(.,t) + N_{\Gamma^*\text{AGFO}} \hat{z}(\xi,t)$$

with

$$M_{\Gamma^*\text{AGFO}} = 0 \text{ and } N_{\Gamma^*\text{AGFO}} = I_{\Gamma^*\text{AGFO}},$$

then we have

$$\hat{x}(\xi,t) = \hat{z}(\xi,t)$$

Now, we can calculate the error of gradient state estimator

$$\hat{e}_{\Gamma^*\text{AGFO}}(\xi,t) = x(\xi,t) - \hat{x}(\xi,t) = x(\xi,t) - x(\xi,t) + x(\xi,t) - \hat{x}(\xi,t)$$

$$= \hat{x}(\xi,t) - \hat{z}(\xi,t) = e(\xi,t) = (A - H_{\Gamma^*\text{AGFO}} C) e(0,t)$$

is asymptotically $\Gamma^*G$-stable with

$$e(0,t) = z_0(\xi) - \hat{z_0}(\xi).$$

Consequently, we get

$$\lim_{t \to \infty} \|x(.,t) - \hat{x}(\xi,t)\|_{(H^{1/2}(\Gamma^*))^n} = \lim_{t \to \infty} \|x(.,t) - \hat{z}(\xi,t)\|_{(H^{1/2}(\Gamma^*))^n} = 0$$

Then, the dynamical system (22) is $\Gamma^*AGFO$-Observer to the system (19).

**Corollary 3.6** From the previous results, we can deduce that**:**

1. Theorem 3.5 gives the sufficient conditions which guarantee the dynamic system (22) is a $\Gamma^*GFO$-observer for the system (19).

2. If a system which is an $\Omega GFO$-observer, then it is $\Gamma^* AGFO$-observer for system (19).

3. If a system is $\Gamma^* GFO$-observer, then it is $\Gamma^{1*} GFO$-observer for every subset $\Gamma^{1*}$ of $\Gamma^*$, but the converse is not true [28]. This is an important case because there are many problems in real world can reconstruction the current state gradient but is not in the usual sense which is presented in the following example.

**Example 3.7.** Consider a two-dimensional system described by the diffusion equation

$$\begin{cases} \frac{\partial x}{\partial t}(\xi_1, \xi_2, t) = \Delta x(\xi_1, \xi_2, t) & Q \\ x(\xi_1, \xi_2, t) = x_0(\xi_1, \xi_2) & \overline{\Omega} \\ \frac{\partial x}{\partial v}(\eta_1, \eta_1, t) = 0 & \Theta \end{cases} \quad (25)$$

where $\Omega = ]0,1[ \times ]0,1[$. The operator $A = \Delta$ generates a strongly continuous semi-group $(S_A(t))_{t \geq 0}$ on the Hilbert space $H^1(\Omega)$ given by

$$S_A(t)x = \sum_{n,m=0}^{\infty} e^{\lambda_{nm} t} \langle x, \varphi_{nm} \rangle_{H^1(\Omega)} \varphi_{nm}$$

With

$$\lambda_{nm} = -(n^2, m^2)\pi^2, \varphi_{nm}(\xi_1, \xi_2) = 2a_{nm} \cos(n\pi\xi_1) \cos(n\pi\xi_2) \text{ and}$$

$$2a_{nm} = (1 - \lambda_{nm})^{-1/2}.$$

Consider the dynamical system

$$\begin{cases} \frac{\partial x}{\partial t}(\xi_1, \xi_2, t) = \Delta x(\xi_1, \xi_2, t) - H_{AGFO} C(x(\xi_1, \xi_2, t) - z(\xi_1, \xi_2, t)) & Q \\ x(\xi_1, \xi_2, 0) = x_0(\xi_1, \xi_2) & \overline{\Omega} \\ \frac{\partial x}{\partial v}(\eta_1, \eta_1, t) = 0 & \Theta \end{cases} \quad (26)$$

Where $H \in L(R^q, Z), Z$ is a Hilbert space and $C: H^1(\overline{\Omega}) \to R^q$ is a linear operator. Consider the boundary sensor $(\Gamma_0, f)$ defined by

$$\Gamma_0 = \{0\} \times ]0,1[ \text{ and } f(\eta_1, \eta_2) = \cos \pi \eta_2$$

Thus, the output function can be written by

$$y(t) = \int_{\Gamma_0} x(\eta_1, \eta_2, t) f(\eta_1, \eta_2) d\eta_1 d\eta_2 \quad (27)$$

If the state $x_0$ is defined in $\Omega$ by

$$x_0(\xi_1, \xi_2) = \cos(\pi \xi_1) \cos(2\pi \xi_2),$$

then the system (20)-(22) is not W$\Omega$G-observable, i.e. the sensor $(\Gamma_0, f)$ is not $\Omega$-strategic and therefore the system (20)-(22) is not $\Omega$AG-detectable [7]. Thus, the dynamical system (21) is not $\Omega$AG-observer for the system (20)-(22) (see [29]). Here, we consider the region $\Gamma = ]0,1[ \times \{0\} \subset \partial \Omega$ (figure 3) and the dynamical system

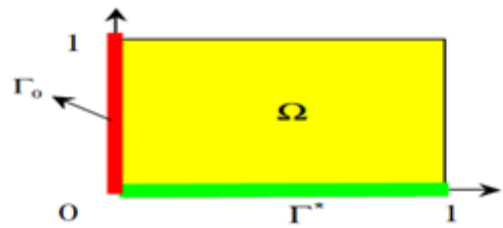

**Fig. 3:** Domain $\Omega$, region $\Gamma^*$ and location $\Gamma_0$ of zone sensor.

$$\begin{cases} \frac{\partial x}{\partial t}(\xi_1,\xi_2,t) = \Delta x(\xi_1,\xi_2,t) \\ \qquad\qquad -H_{\Gamma^*\text{AGFO}\Gamma}C\big(x(\xi_1,\xi_2,t)-z(\xi_1,\xi_2,t)\big) & Q \\ x(\xi_1,\xi_2,0)=x_0(\xi_1,\xi_2) & \bar{\Omega} \\ \frac{\partial x}{\partial v}(\eta_1,\eta_1,t)=0 & \Theta \end{cases} \quad (28)$$

where $H \in L(R^q, H^{1/2}(\Gamma^*))$. In this case, the system (20)-(22) is W$\Gamma^*$G-observable and the sensor $(\Gamma_0, f)$ is $\Gamma^*G$-strategic [26]. Thus, the system (20)-(22) is $\Gamma^*AG$-detectable [13]. Finally the dynamical system (20) is $\Gamma^*AGFO$-observer for the system (20)-(22) [30].

## 4. Application to asymptotic $\Gamma^*AGFO$ -observer in diffusion system

In this section we consider the distributed diffusion systems defined in the domain $\Omega$. Various results related to different types of sensor have been extended. In the case of two-dimensional, we take where $\Omega = ]0,a_1[ \times ]0,a_2[$ and $\Gamma^* = ]0,a_1[ \times \{a_2\}$ is a region of $\partial\Omega$ the boundary of $\Omega$. The eigenfunctions of the dynamic system (16) for Dirichlet boundary conditions are given by

$$\varphi_{nm}(\mu_1,\mu_2) = \left(\frac{4}{a_1 a_2}\right)^{1/2} \cos n\pi \left(\frac{\mu_1}{a_1}\right) \cos m\pi \left(\frac{\mu_2}{a_2}\right) \qquad (29)$$

associated with eigenvalues

$$\lambda_{nm} = -\left(\frac{n^2}{a_1^2}+\frac{m^2}{a_2^2}\right)\pi^2, \ n,m \geq 1 \qquad (30)$$

If we suppose that $a_1^2/a_2^2 \notin Q$, then the multiplicity of the eigenvalues $\lambda_{nm}$ is $r_{nm} = 1$ for every $n,m = \{1,\ldots,J\}$, then one sensor $(D,f)$ may be sufficient for $\Gamma^*AGFO$-observer [24].

### 4.1. Zone case in rectangular domain

This subsection related to give sufficient conditions which is characterized some cases of the $\Gamma^*AGFO$ -observer in the rectangular domain of system (12) with various sensor locations cases.

### 4.1.1. Internal zone sensors case

Consider two dimensional system defined in $\bar{\Omega} = [0,a_1] \times [0,a_2]$ by parabolic equation

$$\begin{cases} \frac{\partial x}{\partial t}(\xi_1,\xi_2,t) = \frac{\partial^2 x}{\partial \xi_1^2}(\xi_1,\xi_2,t) + \frac{\partial^2 x}{\partial \xi_2^2}(\xi_1,\xi_2,t) & Q \\ x(\xi_1,\xi_2,0) = x_0(\xi_1,\xi_2) & \bar{\Omega} \\ x(\eta_1,\eta_2,t)=0 & \Sigma \end{cases} \quad (31)$$

augmented with output function measured by internal or boundary zone sensors

$$y(.,t) = \int_D x(\xi_1,\xi_2,t)f(\xi_1,\xi_2)d\xi_1 d\xi_2 \tag{32}$$

where the zone sensor is the located inside the domain $\Omega$, over the support

$$D = \,]\xi_{1_0} - l_1, \xi_{1_0} + l_1[\,\times\,]\xi_{2_0} - l_2, \xi_{2_0} + l_2[\, \subset \Omega, \text{ and } f \in L^2(D)$$

as in (figure 4).

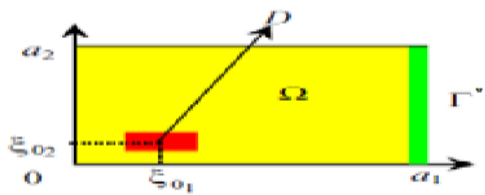

**Fig. 4:** Rectangular domain $\Omega$, region $\Gamma^*$ and location $D$ of internal zone sensor.

In this case the system (31)-(32) have an associated dynamical system is given by the following form

$$\begin{cases} \frac{\partial z}{\partial t}(\xi_1,\xi_2,t) = \frac{\partial^2 z}{\partial \xi_1^2}(\xi_1,\xi_2,t) + \frac{\partial^2 z}{\partial \xi_2^2}(\xi_1,\xi_2,t) \\ \qquad -H_{\Gamma^* G}\big(Cz(\xi_1,\xi_2,t) - y(t)\big) & Q \\ z(\xi_1,\xi_2,0) = z_0(\xi_1,\xi_2) & \bar{\Omega} \\ z(\eta_1,\eta_2,t) = 0 & \Sigma \end{cases} \tag{33}$$

Thus have the following important result.

**Proposition 4.1:** Suppose that $f_1$ is symmetric about $\xi = \xi_{01}$ and $f_2$ is symmetric about $\xi = \xi_{02}$, then the dynamical system (33) is $\Gamma^* AGFO$-observer for systems (31)-(32) if $n\mu_{1_0}/a_1$ and $m\mu_{2_0}/a_2 \notin N$, for every $n,m = \{1,\dots,J\}$.

### 4.1.2. One side boundary zone sensor case

Now, the measurements are given by the following output function

$$y(t) = \int_{\Gamma_0} \frac{\partial x}{\partial v}(\eta_1,\eta_2,t)f(\eta_1,\eta_2)d\eta_1\eta_2, \tag{34}$$

Where $\Gamma_0 \subset \partial\Omega$ is the boundary support of the sensor and $f \in L^2(\Gamma_0)$. In the case, where the support of the sensor $(\Gamma_0, f)$ is one of side as in (figure 5)

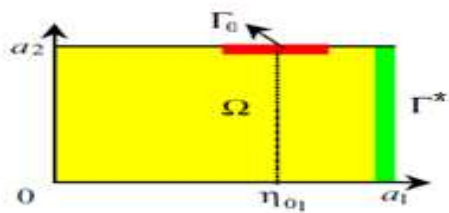

**Fig. 5:** Rectangular domain $\Omega$, region $\Gamma^*$ and location $\Gamma_0$ of boundary zone sensor.

Then, we have the following result.

**Proposition 4.2:** Suppose that $\Gamma^* \subset \partial\Omega$ and $f$ is symmetric with respect to $\eta_1 = \eta_{0_1}$, then the dynamic system (33) is $\Gamma^*AGFO$ -observer for the systems (31)-(32) if $n\eta_{1_0}/a_1 \notin N$, for all $n, \ n = \{1, \ldots, J\}$.

### 4.1.3. Two side boundary zone sensor case

In this case the output function (2) is given by

$$y(t) = \int_{\bar{\Gamma}} \frac{\partial x}{\partial v}(\eta_1, \eta_2, t) f(\eta_1, \eta_2) d\eta_1 \eta_2 \qquad (35)$$

when the support of the sensor is on two sides boundary, i.e.,

$$\bar{\Gamma} = [0, \bar{\eta}_{1_0} + l_1] \times \{0\} \cup \{0\} \times [0, \bar{\eta}_{2_0} + l_2] = \bar{\Gamma}_1 \cup \bar{\Gamma}_2 \subset \partial\Omega,$$

as in (figure 6),

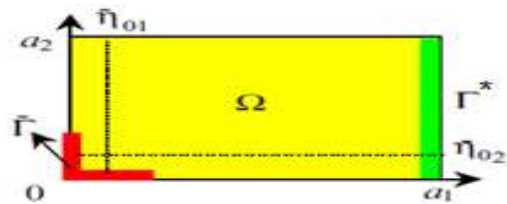

**Fig. 6:** Rectangular domain $\Omega$, region $\Gamma^*$ and location $\bar{\Gamma}^*$ of boundary zone sensor.

then, we obtain the following result.

**Proposition 4.3:** Let $\bar{\Gamma}_1 \cup \bar{\Gamma}_2 \subset \partial\Omega$, the function $f|_{\bar{\Gamma}_1}$ is symmetric with respect to $\eta_1 = \bar{\eta}_{0_1}$, and the function $f|_{\bar{\Gamma}_2}$ is symmetric with respect to $\eta_2 = \bar{\eta}_{0_2}$, then the dynamic system (33) is $\Gamma^*AGFO$ -observer for the systems (31)-(32) if $n\eta_{1_0}/a_1 \notin N$, for all $n, \ n = \{1, \ldots, J\}$.

### 4.2.1. Pointwise sensors case

Consider again the systems (31)-(2) augmented with output function measured by internal pointwise sensors. In this case, the output function is given by the following form

$$y(t) = x(\xi_1, \xi_2, t)\, \delta(\xi_1 - b_1, \xi_2 - b_2) d\xi_1 d\xi_2 \qquad (36)$$

where $b = (b_1, b_2)$ is the location of pointwise sensor in $\Omega$ as defined in (figure 7).

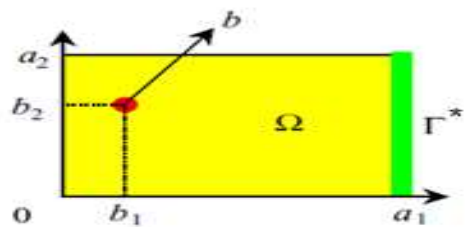

**Fig. 7:** Rectangular domain $\Omega$, region $\Gamma^*$ and location $b$ of internal pointwise sensors.

Then we obtain the following result.

**Proposition 4.4:** Let the sensor is located in $b = (b_1, b_2)$, then, the dynamic system (33) is $\Gamma^*AGFO$ -observer for the systems (31)-(36), if $nb_1/a_1$ and $mb_2/a_2 \notin N$, for every $n, m = \{1, \ldots, J\}$.

### 4.2.2. Filament pointwise case

Suppose that the observation is given by the filament sensor where σ = $Im(\gamma) \subset \Omega$ is symmetric with respect to the line $b = (b_1, b_2)$ as in (figure 8). More precisely, the sensor is line of pointwises inside the domain Ω, then the output function still given by equation (31

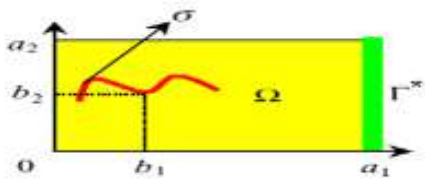

**Fig. 8:** Rectangular domain, region Γ* and locations σ of filament pointwise sensor.

**Proposition 4.5:** Let the sensor is located in $b = (b_1, b_2)$, then, the dynamic system (33) is $\Gamma^* AGFO$ -observer for the systems (31)-(36), if $nb_1/a_1$ and $mb_2/a_2 \notin N$, for every $n, m = \{1, \dots, J\}$.

### 4.2.3 Boundary pointwise case

Suppose that the sensor $(b, \delta_b)$ is located on $b$, where $b = (b_1, b_2) \in \partial\Omega$ with $b = (0, b_2)$ as in (figure 9).

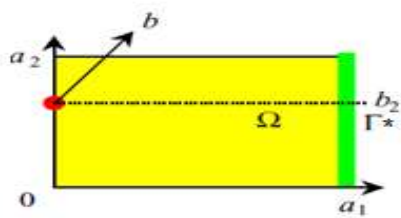

**Fig. 9:** Rectangular domain, region Γ* and location $b$ of boundary pointwise sensor.

In this case, the output function is given by

$$y(.,t) = \int_D x(\eta_1, \eta_2, t)\delta_b(\eta_1, \eta_2 - b_2)d\eta_1\eta_2 (y(.,t)$$

$$= \int_{\partial\Omega} \frac{\partial x}{\partial v}(\eta_1, \eta_2, t)f(\eta_1, \eta_2)\, d\eta_1 d\eta_2 \qquad (37)$$

where $b = (0, b_2)$. Thus, we obtain the following result.

**Proposition 4.6:** Let the sensor is located in $b = (0, b_2)$, then, the dynamic system (33) is $\Gamma^* AGFO$-observer for the systems (31)-(37), if $mb_2/a_2 \notin N$ for all $m = \{1, \dots, J\}$.

Now, from the previous results in this paper, we can deduce the following results.

**Remark 4.7**: We can extend these results to the case of two dimensional systems with circular domain in different sensor structures as in [13, 16].

**Remark 4.8**: We can extend the above results of the two dimensional systems (31)-(2) to case of one dimensional systems case with $\Omega = ]0, a[$ as in ref.s [8-17, 25].

**Remark 4.9**: We know that the previous results have been developed with Dirichlet boundary conditions, then we can extend with Neumann or mixed boundary conditions as in [10, 31].

## 5. Conclusion

The concept have been studied in this paper is related to the $\Gamma^*AGFO$-observer in connection with sensors structure for a class of distributed parameter systems. More precisely, we have been given a sufficient condition for the existing an $\Gamma^*AGFO$-observer which allows to estimate the gradient state in a subregion $\Gamma^*$. For future work, one can be extension these result to the problem of regional boundary asymptotic gradient reduced order observer in connection with the sensors structures.

**Acknowledgments.** Our thanks in advance to the editors and experts for their efforts.